\documentclass[a4paper]{amsart}
\usepackage{amssymb}
\usepackage{ifthen}
\usepackage{graphicx}

\newtheorem{thm}{Theorem}
\newtheorem{cor}{Corollary}
\newtheorem{lem}{Lemma}

\newtheorem{rem}{Remark}

\newtheorem{conj}{Conjecture}
\theoremstyle{definition}
\newtheorem{example}[equation]{Example}%[section]
\newtheorem{prob}[equation]{Problem}

\newcommand{\ID}{{\mathbb D}}

\newcommand{\D}{{\mathbb D}}
\def\be{\begin{equation}}
\def\ee{\end{equation}}

\newcommand{\bee}{\begin{enumerate}}
\newcommand{\eee}{\end{enumerate}}

\newcommand{\blem}{\begin{lem}}
\newcommand{\elem}{\end{lem}}
\newcommand{\bthm}{\begin{thm}}
\newcommand{\ethm}{\end{thm}}
\newcommand{\bcor}{\begin{cor}}
\newcommand{\ecor}{\end{cor}}
\newcommand{\beg}{\begin{example}}
\newcommand{\eeg}{\end{example}}
\newcommand{\begs}{\begin{examples}}
\newcommand{\eegs}{\end{examples}}
\newcommand{\bdefe}{\begin{defin}}
\newcommand{\edefe}{\end{defin}}
\newcommand{\bprob}{\begin{prob}}
\newcommand{\eprob}{\end{prob}}
\newcommand{\bei}{\begin{itemize}}
\newcommand{\eei}{\end{itemize}}

\newcommand{\bcon}{\begin{conj}}
\newcommand{\econ}{\end{conj}}
\newcommand{\bcons}{\begin{conjs}}
\newcommand{\econs}{\end{conjs}}
\newcommand{\bprop}{\begin{propo}}
\newcommand{\eprop}{\end{propo}}
\newcommand{\br}{\begin{rem}}
\newcommand{\er}{\end{rem}}
\newcommand{\brs}{\begin{rems}}
\newcommand{\ers}{\end{rems}}
\newcommand{\bo}{\begin{obser}}
\newcommand{\eo}{\end{obser}}
\newcommand{\bos}{\begin{obsers}}
\newcommand{\eos}{\end{obsers}}
\newcommand{\bpf}{\begin{pf}}
\newcommand{\epf}{\end{pf}}
\newcommand{\ba}{\begin{array}}
\newcommand{\ea}{\end{array}}
\newcommand{\beq}{\begin{eqnarray}}
\newcommand{\beqq}{\begin{eqnarray*}}
\newcommand{\eeq}{\end{eqnarray}}
\newcommand{\eeqq}{\end{eqnarray*}}

\begin{document}
\bibliographystyle{amsplain}

\title[Sharp bounds of logarithmic coefficients ]{Sharp bounds of logarithmic coefficients for a class
of univalent functions}

\author[M. Obradovi\'{c}]{Milutin Obradovi\'{c}}
\address{Department of Mathematics,
Faculty of Civil Engineering, University of Belgrade,
Bulevar Kralja Aleksandra 73, 11000, Belgrade, Serbia}
\email{obrad@grf.bg.ac.rs}

\author[N. Tuneski]{Nikola Tuneski}
\address{Department of Mathematics and Informatics, Faculty of Mechanical Engineering, Ss. Cyril and
Methodius
University in Skopje, Karpo\v{s} II b.b., 1000 Skopje, Republic of North Macedonia.}
\email{nikola.tuneski@mf.edu.mk}

%\date{\today  %June. 30, 09
%;  File: }

%\begin{abstract}

%\end{abstract}

%\thanks{The work of the first author was supported by MNZZS Grant, No. ON174017, Serbia. The research of
%the second
%author was supported by National Board for Higher
%Mathematics, India.}

\begin{abstract}
Let $\mathcal{U(\alpha, \lambda)}$, $0<\alpha <1$, $0 < \lambda <1$ be the class of
functions $f(z)=z+a_{2}z^{2}+a_{3}z^{3}+\cdots$ satisfying
$$\left|\left(\frac{z}{f(z)}\right)^{1+\alpha}f'(z)-1\right|<\lambda$$
in the unit disc $\ID$. For $f\in \mathcal{U(\alpha, \lambda)}$ we give sharp bounds of its initial
logarithmic coefficients $\gamma_{1},\,\gamma_{2},\,\gamma_{3}.$
\end{abstract}
\maketitle
%\pagestyle{myheadings}
%\markboth{M. Obradovi\'{c} and  S. Ponnusamy }{}
%\cc
%\section{Introduction}

\medskip

\section{Introduction and definitions}

\medskip

Let $\mathcal{A}$ be the class of functions $f$ which are analytic  in the open unit disc $\D=\{z:|z|<1\}$
of the form
\be\label{eq1}
f(z)=z+a_2z^2+a_3z^3+\cdots,
\ee
and let $\mathcal{S}$ be the subclass of $\mathcal{A}$ consisting of functions that are univalent in $\ID$.

\medskip

For a function $f\in \mathcal{S}$ we define its logarithmic coefficients, $\gamma_n$, $n=1,2,\ldots$,  by
\be\label{eq2}
\log\frac{f(z)}{z}=2\sum_{n=1}^\infty \gamma_n z^n.
\ee
Relatively little exact information is known about those coefficients. The natural conjecture
$|\gamma_n|\le1/n$, inspired by the Koebe function (whose logarithmic coefficients are $1/n$) is false
even in order of magnitude (see Duren \cite{duren}).
For the class $\mathcal{S}$ the sharp estimates of single logarithmic coefficients  are known only for
$\gamma_1$ and $\gamma_2$, namely,
\[|\gamma_1|\le1\quad\mbox{and}\quad |\gamma_2|\le \frac12+\frac1e=0.635\ldots,\]
and are unknown for $n\ge3$. The best known estimate $|\gamma_3|\leq 0.55661\ldots$ was given by the
authors (see \cite{MONT}).
For the subclasses of univalent functions  the situation is not a great deal better. Only the estimates of
the initial logarithmic coefficients are available. For details see \cite{cho}.

\medskip

In the paper \cite{Ob1}   the class $\mathcal{U}(\alpha, \lambda)$  ($0<\alpha <1$, $0 < \lambda < 1 $)
of functions $f\in \mathcal{A}$  was introduced by the condition
\be\label{eq3}
\left|\left(\frac{z}{f(z)}\right)^{1+\alpha}f'(z)-1\right|<\lambda, \quad z\in\ID.
\ee
There is  shown that functions from $\mathcal{U}(\alpha, \lambda)$ are starlike, i.e., belong to the class $\mathcal{S}^{\star}$ of functions that map the unit disk onto a starlike domain, if
\be\label{eq4}
0<\lambda \leq \frac{1-\alpha}{\sqrt{(1-\alpha)^{2}+\alpha^{2}}} \equiv \lambda_{\star}.
\ee
In the limiting cases when $\lambda=1,$ and either $\alpha=0$ or $\alpha=1,$ functions in the classes $\mathcal{U}(0,1)$ and $\mathcal{U}(1,1)$ satisfy
$$ \left|\frac{zf'(z)}{f(z)}-1\right|<1, \quad \text{and} \quad \left|\left(\frac{z} {f(z)}\right)^2f'(z)-1\right|<1,$$
respectively. The former is a subclass of $\mathcal{S}^{\star}$ since the analytical characterisation of starlike functions is $\operatorname{Re}\frac{zf'(z)}{f(z)}>0$ ($z\in\ID$), while functions in the latter class are univalent (see \cite{OP-2011, MONT-2019}).

\medskip

In this paper we consider estimates of three initial logarithmic coefficients for the class
$\mathcal{U}(\alpha, \lambda)$, where  $0<\alpha <1, \, 0 < \lambda\leq\lambda_{\star}$ and $\lambda_{\star}$ is defined by \eqref{eq4}.

\medskip

For our consideration we need the next lemma.

\begin{lem}\label{lem1}\cite{Ob2} Let $f\in \mathcal{U}(\alpha, \lambda),\, 0<\alpha <1, \, 0 <\lambda< 1.$ Then there exists a function  $\omega$, analytic in $\ID$, such that $\omega(0)=0$, $| \omega(z)|<1$ for all $z\in\ID$, and
\be\label{eq5}
\left[\frac{z}{f(z)}\right]^{\alpha}= 1-\alpha \lambda z^{\alpha}\int_{0}^{z}\frac{\omega(t)}{t^{\alpha+1}}dt.
\ee
\end{lem}

\medskip

By $\Omega $ we denote the class of analytic functions in $\ID$:
\be\label{eq6}
\omega(z)=c_{1}z+c_{2}z^2+c_{3}z^3+\cdots,
\ee
with $\omega(0)=0$, and $| \omega(z)|<1$ for all $z\in\ID$.

\medskip

In their paper \cite{Pro} Prokhorov and Szynal  obtained sharp estimates on the functional $$\Psi(\omega)=|c_{3}+\mu c_{1}c_{2}+\nu c_{1}^{3}|$$
within the class of all $\omega \in \Omega.$
For our application we need only a part of those results.

\medskip

\begin{lem}\label{lem2}\cite{Pro}
Let $\omega(z)=c_{1}z+c_{2}z^2+c_{3}z^3+\cdots \in \Omega.$ For  $\mu$ and $\nu$ real numbers, let
$$\Psi(\omega)=\left|c_{3} + \mu c_{1}c_{2} + \nu c_{1}^3 \right|,$$
and
\beqq
D_{1}& =& \left\{(\mu,\nu): |\mu| \leq \frac{1}{2},  |\nu| \leq 1 \right\},\\
D_{2} &=& \left\{(\mu,\nu): \frac{1}{2} \leq |\mu| \leq 2, \frac{4}{27}(|\mu|+1)^3 - (|\mu|+1) \leq \nu \leq 1 \right\},\\[2mm]
D_{3}& =& \left\{(\mu,\nu): |\mu| \leq 2, |\nu| \geq 1 \right\}.
\eeqq
Then, the sharp estimate $\Psi(\omega) \leq \Phi(\mu,\nu)$ holds, where
\[\Phi(\mu,\nu) = \left\{
  \begin{array}{cc}
    1, & \hbox{$(\mu,\nu) \in D_{1} \cup D_{2} \cup \{(2,1)\}$;} \\
    |\nu|, & \hbox{$(\mu,\nu) \in D_{3} $.}
  \end{array}
\right.\]
\end{lem}

\medskip

\section{Main results}

\bthm\label{23-th 1}
Let $f(z)=z+a_{2}z^{2}+a_{3}z^{3}+ \cdots$ belongs to the class $\mathcal{U}(\alpha, \lambda)$ and
$\lambda_{\star}$ is defined by \eqref{eq4}. Then the following results are best possible.
\begin{itemize}
  \item[($i$)]   $|\gamma_{1}|\leq \frac{\lambda }{2(1-\alpha)}$ when  $0<\lambda\leq \lambda_{\star}$ and $0<\alpha<1.$
  \item[($ii$)]  Let $\lambda _{1}= \frac{2(1-\alpha)^{2}}{\alpha(2-\alpha)}$ and let $\alpha_{1}=0.4825\ldots$ be the unique real root of the equation
$$7\alpha^{4}-20\alpha^{3}+24\alpha^{2}-16\alpha+4=0$$
on the interval $(0,1)$. Then
$$|\gamma_{2}|\leq \frac{\lambda }{2(2-\alpha)} \quad \text{if} \quad 0<\lambda \leq \begin{cases} \lambda_{1},\,\alpha\in[\alpha_{1},1),\\
\lambda_{\star},\, \alpha\in(0,\alpha_{1}],\end{cases}$$
and
$$|\gamma_{2}|\leq \frac{\alpha\lambda^{2}}{4(1-\alpha)^{2}}  \quad \text{if} \quad
\lambda _{1}\leq \lambda \leq \lambda_{\star},\,\alpha\in[\alpha_{1}, 1).$$
\item[($iii$)]  Let $ \lambda_{1/2}=\frac{(1-\alpha)(2-\alpha)}{2\alpha(3-\alpha)},$
$\lambda_{\nu}=\sqrt{\frac{3(1-\alpha)^{3}}{\alpha^{2}(3-\alpha)}}$
and $\alpha_{1/2}=0.2512\ldots$ and $\alpha_{\nu}=0.5337\ldots$
are the unique roots of equations
$$ 4-12\alpha-19\alpha^{2}+14\alpha^{3}-2\alpha^{4}=0 $$
and
$$3-9\alpha+9\alpha^{2}-5\alpha^{3}=0,$$
on the interval $(0,1)$, respectively. Then
$$|\gamma_{3}|\leq \frac{\lambda}{2(3-\alpha)} \quad \text{if}\quad
0<\lambda \leq \left\{\begin{array}{cc} \lambda_{\star},& \alpha\in(0,\alpha_{1/2}],\\
\lambda_{1/2},& \alpha\in[\alpha_{1/2},\alpha_{2}], \\
\lambda_{\nu}, & \alpha\in[\alpha_{2},1),\end{array}\right.$$
where $\alpha_{2}=0.9555\ldots$ is the unique real root of equation
$11\alpha^{2}-44\alpha+32=0$ on $(0,1)$. Also,
$$|\gamma_{3}|\leq \frac{\alpha^{2}\lambda^{3}}{6(1-\alpha)^{3}}  \quad \text{if}\quad
\lambda_{\nu}\leq \lambda \leq \lambda_{\star},\,\alpha\in[\alpha_{\nu},1).$$
\end{itemize}
%All these results are the best possible.
\ethm

\begin{proof}
Let $f\in \mathcal{U}(\alpha, \lambda)$ and  $\omega\in \Omega$ are given by \eqref{eq1} and
\eqref{eq6}, respectively. Then, from \eqref{eq5}, upon integration, we have
$$ \left[\frac{z}{f(z)}\right]^{\alpha}= 1-\alpha \lambda \sum_{n=1}^{\infty}\frac{c_{n}}{n-\alpha}z^{n},$$
that is,
\be\label{eq7}
\frac{f(z)}{z}=\left(1-\alpha \lambda \sum_{n=1}^{\infty}\frac{c_{n}}{n-\alpha}z^{n}\right)^{-\frac{1}{\alpha}}
\ee
(the principal value is used here).
Further, from \eqref{eq7}, having in mind that
$$(1-\alpha z)^{-1/\alpha} = 1+z+\frac{1+\alpha}{2} z^2+\frac{(1+\alpha) (1+2\alpha)}{6}  z^3+\cdots,$$
after some calculations, we obtained
\begin{equation*}
  \begin{split}
\sum_{n=1}^{\infty}a_{n+1}z^{n}&=\sum_{n=1}^{\infty}\frac{\lambda c_{n}}{n-\alpha}z^{n}+
\frac{1+\alpha}{2}\left(\sum_{n=1}^{\infty}\frac{\lambda c_{n}}{n-\alpha}z^{n}\right)^{2}\\
&+ \frac{(1+\alpha)(1+2\alpha)}{6}\left(\sum_{n=1}^{\infty}\frac{\lambda c_{n}}{n-\alpha}z^{n}\right)^{3}+\cdots.
  \end{split}
\end{equation*}
By comparing the coefficients we receive
\be\label{eq8}
\begin{split}
a_{2}&=\frac{\lambda }{1-\alpha}c_{1}, \\[2pt]
a_{3}&=\frac{\lambda }{2-\alpha}c_{2}+\frac{(1+\alpha)\lambda^{2}}{2(1-\alpha)^{2}}c_{1}^{2},\\[2pt]
a_{4}&=\frac{\lambda}{3-\alpha}c_{3}+\frac{(1+\alpha)\lambda^{2}}{(1-\alpha)(2-\alpha)}c_{1}c_{2}
+\frac{(1+\alpha)(1+2\alpha)\lambda^{3}}{6(1-\alpha)^{3}}c_{1}^{3}.
\end{split}
\ee
On the other hand, by comparing the coefficients  in the relation \eqref{eq2}, for the logarithmic coefficients we obtain
\be\label{eq9}
\gamma_{1}=\frac{1}{2}a_{2},\quad \gamma_{2}=\frac{1}{4}(2a_{3}-a_{2}^{2}),\quad
\gamma_{3} =\frac{1}{2}(a_4-a_2a_3+\frac{1}{3}a_{2}^{3}).
\ee
Using the relations \eqref{eq8} and \eqref{eq9}, after some calculations, we have
\be\label{eq10}
\begin{split}
\gamma_{1}&=\frac{\lambda }{2(1-\alpha)}c_{1}, \\[2pt]
\gamma_{2}&=\frac{1}{4}\left[\frac{2\lambda }{2-\alpha}c_{2}+\frac{\alpha\lambda^{2}}{(1-\alpha)^{2}}c_{1}^{2}\right],\\[2pt]
\gamma_{3}&=\frac{\lambda}{2(3-\alpha)}\left(c_{3}+\mu c_{1}c_{2}
+\nu c_{1}^{3}\right),
\end{split}
\ee
where
\be\label{eq11}
\mu=\frac{\alpha(3-\alpha)\lambda}{(1-\alpha)(2-\alpha)}\quad \text{and}\quad \nu=\frac{\alpha^{2}(3-\alpha)\lambda^{2}}{3(1-\alpha)^{3}}.
\ee

\medskip

Since logarithmic coefficients are defined for univalent functions, in order to guarantee univalence of $f$ in all cases we need $ 0<\lambda\leq \lambda_{\star}$, where $\lambda_{\star}$ is defined in \eqref{eq4}.

\medskip

\begin{itemize}
  \item[($i$)]
 From \eqref{eq10} we have $|\gamma_{1}|\leq \frac{\lambda }{2(1-\alpha)}$, where
$0<\lambda\leq \lambda_{\star}$ and $0<\alpha<1.$ The result is the best possible as the function $f_{1}$
defined by
$$ f_{1}(z)=z\left( 1-\frac{\alpha \lambda}{1-\alpha}z\right)^{-1/\alpha}=z+\frac{ \lambda}{1-\alpha}z^{2}+\ldots$$
shows.
\medskip

  \item[($ii$)] Using the inequalities $|c_{1}|\leq 1,\,|c_{2}|\leq 1-|c_{1}|^{2}$ for $\omega\in \Omega$ and
\eqref{eq10}, we have
\beqq
\begin{split}
|\gamma_{2}|&\leq \frac{1}{4}\left[\frac{2\lambda }{2-\alpha}|c_{2}|+\frac{\alpha\lambda^{2}}{(1-\alpha)^{2}}|c_{1}|^{2}\right]\\[2mm]
&\leq\frac{1}{4}\left[\frac{2\lambda }{2-\alpha}(1-|c_{1}|^{2})+\frac{\alpha\lambda^{2}}{(1-\alpha)^{2}}|c_{1}|^{2}\right]\\[2mm]
&\leq\frac{1}{4}\left[\frac{2\lambda }{2-\alpha}+\left(\frac{\alpha\lambda^{2}}{(1-\alpha)^{2}}-\frac{2\lambda }{2-\alpha}\right)|c_{1}|^{2}\right]\equiv H_1(|c_1|).
\end{split}
\eeqq
 If $\frac{\alpha\lambda^{2}}{(1-\alpha)^{2}}-\frac{2\lambda }{2-\alpha}\leq 0$, or equivalently,
$$\lambda \leq\frac{2(1-\alpha)^{2}}{\alpha(2-\alpha)}\equiv \lambda _{1},$$ then
$|\gamma_{2}|\leq H_1(0) = \frac{\lambda }{2(2-\alpha)}$. It is also necessary that
$$\lambda\leq \lambda_{\star}=\frac{1-\alpha}{\sqrt{(1-\alpha)^{2}+\alpha^{2}}}.$$
The last inequality will hold if $\lambda_1\le\lambda_\star$, or equivalently, if
$$7\alpha^{4}-20\alpha^{3}+24\alpha^{2}-16\alpha+4\leq0,$$
i.e., if  $\alpha\in[\alpha_{1}, 1)$, where $\alpha_{1}=0.4825\ldots$ is the unique real root of equation
$$7\alpha^{4}-20\alpha^{3}+24\alpha^{2}-16\alpha+4=0$$
on the interval $(0,1)$. If $\alpha\in(0,\alpha_{1}]$, then $\lambda _{1}\geq \lambda_{\star}$ and we have that
$0<\lambda\leq \lambda_{\star}$ will imply the same result.

\medskip

Finally, if $\alpha\in[\alpha_{1}, 1)$, i.e., $\lambda_1\le\lambda_\star$, and $\lambda _{1}\leq \lambda \leq \lambda_{\star}$, then, from the previous consideration we obtain
$$|\gamma_{2}|\leq H_1(1) = \frac{\alpha\lambda^{2}}{4(1-\alpha)^{2}} .$$

\medskip

Those results are the best possible as the functions given by \eqref{eq7} for $c_{2}=1$ $ (c_{1}=c_{3}=\cdots=0)$ or for $c_{1}=1$ $ (c_{2}=c_{3}=\cdots=0)$, show.
\medskip

  \item[($iii$)] From \eqref{eq10}  we have
\be\label{eq12}
|\gamma_{3}|\leq \frac{\lambda}{2(3-\alpha)}\left|c_{3}+\mu c_{1}c_{2}
+\nu c_{1}^{3}\right|=\frac{\lambda}{2(3-\alpha)}\Psi(\omega),
\ee
where   $\mu $ and $\nu$ are given by \eqref{eq11}.

\medskip

Next, we want to apply the results of Lemma \ref{lem2}, and for that we need to distinguish the cases in the definitions of the sets $D_1$, $D_2$, and $D_3$.

\medskip

First, we note that $\mu $ and $\nu$ are both positive.

\medskip

Further, $ \mu=\frac{\alpha(3-\alpha)\lambda}{(1-\alpha)(2-\alpha)}\leq \frac{1}{2}$ is equivalent to
$$ 0<\lambda\leq \frac{(1-\alpha)(2-\alpha)}{2\alpha(3-\alpha)}\equiv\lambda_{1/2}.$$
It is necessary that $\lambda\leq \lambda_{\star}$, where $\lambda_{\star}$ is defined by \eqref{eq4}. After some calculations,  $\lambda_{1/2}\leq \lambda_{\star}$ is equivalent to
 $$ 4-12\alpha-19\alpha^{2}+14\alpha^{3}-2\alpha^{4}\leq0 ,$$
i.e., to $\alpha\in[\alpha_{1/2},1)$, where $\alpha_{1/2}=0.2512\ldots$ is the unique real root of the equation
$$ 4-12\alpha-19\alpha^{2}+14\alpha^{3}-2\alpha^{4}=0 $$
on the interval $(0,1)$. In that sense we have
\be\label{eq13}
0<\mu\leq \frac{1}{2}\quad \Leftrightarrow \quad \lambda \leq \left\{\begin{array}{cc} \lambda_{1/2}, &\alpha\in[\alpha_{1/2},1),\\
\lambda_{\star}, &\alpha\in(0,\alpha_{1/2}]. \end{array}\right.
\ee
On the other hand, by \eqref{eq11}, $\nu=\frac{\alpha^{2}(3-\alpha)\lambda^{2}}{3(1-\alpha)^{3}}\leq1$ is equivalent to
$$ 0<\lambda\leq \sqrt{\frac{3(1-\alpha)^{3}}{\alpha^{2}(3-\alpha)}}\equiv \lambda_{\nu}.$$
It is again necessary that $\lambda\leq \lambda_{\star}$.

\medskip

Next, $\lambda_\nu\leq \lambda_{\star}$ after some calculations is
equivalent to
$$3-9\alpha+9\alpha^{2}-5\alpha^{3}\leq0,$$ which is true when $\alpha\in[\alpha_{\nu},1)$, where
$\alpha_{\nu}=0.5337\ldots$ is the unique real root of equation
$$3-9\alpha+9\alpha^{2}-5\alpha^{3}=0.$$
It means that
\be\label{eq13}
0<\nu\leq 1 \quad \Leftrightarrow \quad \lambda \leq \left\{\begin{array}{cc} \lambda_{\nu},& \alpha\in[\alpha_{\nu},1),\\
\lambda_{\star}, & \alpha\in(0,\alpha_{\nu}]. \end{array}\right.
\ee
Also, $\lambda_{1/2}\leq \lambda_{\nu}$ is equivalent to $11\alpha^{2}-44\alpha+32\geq0$, i.e., to
$\alpha\in[\alpha_{2},1),$ where $\alpha_{2}=0.9555\ldots$ is the unique real root of equation
$$11\alpha^{2}-44\alpha+32=0$$
on the interval $(0,1)$.

\medskip

Using all those previous facts, we can conclude that if
$$
0<\lambda \leq \left\{\begin{array}{cc} \lambda_{\star},& \alpha\in(0,\alpha_{1/2}],\\
\lambda_{1/2},& \alpha\in[\alpha_{1/2},\alpha_{2}], \\
\lambda_{\nu},& \alpha\in[\alpha_{2},1),\end{array}\right.
$$
then $0<\mu\leq\frac{1}{2}$ and $0<\nu\leq1$. By Lemma \ref{lem2} (case $D_{1}$) it means that $\Psi(\omega)\leq1$
and so, by \eqref{eq12}:
$$|\gamma_{3}|\leq \frac{\lambda}{2(3-\alpha)}.$$
The result is best possible as the function obtained for  $c_{3}=1$ $ (c_{1}=c_{2}=c_{4}=\cdots=0)$ in \eqref{eq7} shows.

\medskip

If $\lambda_{1/2}\leq\lambda\leq\lambda_{\nu}$ $\alpha_{\nu}\leq\alpha \leq\alpha_{2}$, then $0<\nu\leq1$ and
\begin{equation*}
\begin{split}
\frac{1}{2}\leq\mu&=\frac{\alpha(3-\alpha)\lambda}{(1-\alpha)(2-\alpha)}
\leq\frac{\alpha(3-\alpha)\lambda_{\nu}}{(1-\alpha)(2-\alpha)}\\
&=\frac{\sqrt{3(1-\alpha)(3-\alpha)}}{2-\alpha}\leq 1.2667\ldots.
\end{split}
\end{equation*}
The last is obtained for $\alpha=\alpha_{\nu}=0.5337\ldots$ since  $\frac{\sqrt{3(1-\alpha)(3-\alpha)}}{2-\alpha}$ is a decreasing function on $(\alpha_{\nu},\alpha_{2})$.

\medskip

For the study of the set $D_2$, we note that the function
$$\phi(\mu)\equiv\frac{4}{27}(1+\mu)^{3}-(1+\mu)$$
is an increasing function for $\frac{1}{2}\leq\mu\leq2$, and
$$\phi(\mu)\leq\phi(1.2667\ldots)=-0,541\ldots<0<\nu\leq1.$$
That implies $\Psi(\omega)\leq1$ (by Lemma \ref{lem2}, case $D_{2}$), and follows the same sharp estimate $|\gamma_{3}|\leq \frac{\lambda}{2(3-\alpha)}$ as in previous case.

\medskip

Finally, since for all $0<\lambda\leq\lambda_{\star}$ we have $0<\mu\leq2$ (easy to check) and if
$\lambda_{\nu}\leq \lambda \leq \lambda_{\star},\,\alpha\in[\alpha_{\nu},1)$, then by Lemma \ref{lem2} (case $D_{3}$):
$\Psi(\omega)\leq \nu $, which by \eqref{eq12} implies
$$|\gamma_{3}|\leq \frac{\lambda}{2(3-\alpha)}\frac{\alpha^{2}(3-\alpha)\lambda^{2}}{3(1-\alpha)^{3}}
=\frac{\alpha^{2}\lambda^{3}}{6(1-\alpha)^{3}}.$$
The result is the best possible as the function given by \eqref{eq7} and $c_{1}=1$ $ (c_{2}=c_{3}=\cdots=0)$ shows.
\end{itemize}
\end{proof}

\end{document}